\documentclass{article}[10pt]

\usepackage{latexsym}
\usepackage{amssymb}
\usepackage{amsmath,amsfonts}
\usepackage{graphicx}
\usepackage{tikz}
\usetikzlibrary{arrows}
\usepackage{multicol}

\newtheorem{theorem}{Theorem}[section]
\newtheorem{prop}[theorem]{Proposition}
\newtheorem{lemma}{Lemma}[section]
\newtheorem{corollary}{Corollary}[section]
\newtheorem{example}{Example}[section]
\newtheorem{definition}{Definition}[section]

\newtheorem{conjecture}[theorem]{Conjecture}

\newtheorem{question}[theorem]{Question}

\def \dd{\hfill \baseb \vskip .5cm}
\def \d{{\noindent \it Proof. } }

\def \Gn{{{\cal G}_n}}

\def \Z{{{\cal Z}}}

\def \Dn{{{\cal D}_n}}

\def \zk{{\mathbb Z}_k}

\def\ex{\begin{example}}
\def\eex{\end{example}}
\def\exx{\end{example}}
\def\t{\begin{theorem}}
\def\tt{\end{theorem}}
\def\D{\begin{definition}}
\def\DD{\end{definition}}
\def\l{\begin{lemma}}
\def\ll{\end{lemma}}
\def\c{\begin{corollary}}
\def\cc{\end{corollary}}
\def\cj{\begin{conjecture}}
\def\cjj{\end{conjecture}}
\def\e{\begin{equation}}
\def\ee{\end{equation}}
\def\p{\begin{prop}}
\def\pp{\end{prop}}
\def\q{\begin{question}}
\def\qq{\end{question}}

\def \Dn{{\mathcal D}_n}

\def \Gn{{\mathcal G}_n}

\def \Tn{{\mathcal T}_n}

\def\A{{\cal A}}
\def\Z{{\cal Z}}
\def\++{\boxplus}

\def\>{\longrightarrow}
\def\<{\longleftarrow}
\newcommand{\baseb}{\hfill \rule{2mm}{2mm}}

\begin{document}

\baselineskip 15pt

\title{{\bf Cordial Digraphs}\thanks{AMS Classification number: 05C20, 05C38, 05C78
} }
\author{LeRoy B. Beasley\\Clock Tower Plaza,  Ste 317\\ 550 North Main St, Box C3 \\ Logan, Utah 84321, U.S.A\\ leroy\_beas@aol.com}
\date{}

\maketitle
\thispagestyle{empty}

\begin{abstract}
A $(0,1)$-labeling of a set is said to be {\em friendly} if the number of  elements of the set  labeled 0 and the number labeled 1 differ by at most 1.  Let $g$ be a labeling of the edge set of a graph that is induced by a labeling $f$ of the vertex set.  If both $g$ and $f$ are friendly then $g$ is said to be a {\em cordial} labeling of the graph.  We extend this concept to directed graphs and investigate the cordiality of directed graphs.  We show that all directed paths and all directed cycles are cordial.  We also discuss the cordiality of oriented trees and other digraphs.
\end{abstract}

\section{Introduction}
Let $\Gn$ denote the set of all simple undirected graphs on the vertex set $V=\{v_1,v_2,\dots,v_n\}$.  For $G\in\Gn$, let $E(G)$ denote the edge set of $G$.  In the article ``Cordial graphs: A Weaker Version of Graceful and Harmonious Graphs'', Ibrahim Cahit introduced cordial graphs.  Relative definitions are found in the next section.  In this article we are interested in directed graphs and so, the orientation of an arc is essential in any investigation.  While the induced labelings in undirected graphs are symmetric, we shall want the arc labelings of our directed graphs to  be asymmetric.  The study of cordial graphs has been to identify classes of graphs that have a friendly labeling that induces a cordial labeling of the edges.  We shall proceed in that fashion for directed graphs.

\section{Preliminaries}

In this article, we will be concerned mainly with digraphs,  We let $\Dn$ denote the set of all simple directed graphs on the vertex set $V = \{v_1,v_2,\dots,v_n\}$.  Further, we shall want no arcs directed both ways between any pair of vertices, that is our graphs will be digon-free.  Thus we shall let $\Tn$ denote the set of all subgraphs of a tournament digraph.  Let $D\in\Tn$, $D=(V,A)$ where $A$ is the arc set of $D$.  Then $D$  has no loops, no digons, and no multiple arcs.  An arc in $D$ directed from vertex $u$ to vertex $v$ will be denoted $\overrightarrow{uv}, \overleftarrow{vu}$ or by the ordered pair $(u,v)$. We also let $\Gn$ denote the set of all simple undirected graphs on the vertex set $V = \{v_1,v_2,\dots,v_n\}$.  So all members of $\Tn$ are orientations of graphs in $\Gn$.  

 A $(0,1)$-labeling of a finite set, ${\cal Z}$,  is a mapping $f:{\cal Z}\to \{0,1\}$ and is  said to be \underline{\em friendly} if approximately one half of the members of ${\cal Z}$ are labeled 0 and the others are labeled 1, that is $-1\leq |f^{-1}(0)|-|f^{-1}(1)| \leq 1$ where $|{\cal X}|$ denotes the cardinality of the set ${\cal X}$.  If there are an even number of elements in ${\cal Z}$ , then a friendly labeling has the number labeled 1 equal to the number labeled 0.  More generally we define an $\A$-friendly labeling:
\D Let $\A$ be a finite set.  An $\A$-labeling of the set \Z\, is a  mapping $f: \Z  \to \A$.   The labeling $f$ is \underline{\em $\A$-friendly} if $-1\leq |f^{-1}(i)|- |f^{-1}(j)|\leq 1$ for any $i,j\in\A$.  
\DD 
Note that a labeling is \A-friendly means that the labels are as nearly evenly distributed as possible.  If the cardinality of \Z \, is a multiple of the cardinality of \A\, then $ |f^{-1}(i)|= |f^{-1}(j)|$  for any $i,j\in\A$.  If the set \A\, is obvious from the context we just say that $f$ is friendly.

Let $G=(V,E)$ be an undirected graph with vertex set $V$ and edge set $E$, and let $f$ be a friendly (0,1)-labeling of the vertex set $V$.   Given this friendly vertex labeling $f$, an induced $(0,1)$-labeling of the edge set  is a mapping $g:E\to \{0,1\}$ where for an edge $uv,\, g(uv)= \hat{g}(f(u),f(v))$ for some $\hat{g}:\{0,1\}\times\{0,1\}\to \{0,1\}$ and is said to be cordial if $g$ is also friendly, that is  $-1\leq |g^{-1}(0)|-|g^{-1}(1)| \leq 1$.  A graph, $G$,  is called {\em cordial} if there exists a induced cordial labeling of the edge set of $G$.  The induced labeling $g$ is usually $g(u,v)=\hat{g}(f(u),f(v))=|f(v)-f(u)|$ \cite{C},  $g(u,v)=\hat{g}(f(u),f(v))=f(v)+f(u)$ (in ${\mathbb Z}_2$) \cite{H},  or $g(u,v)=\hat{g}(f(u),f(v))=f(v)f(u)$ (product cordiality) \cite{S}.  

In\cite{H}, Hovey introduced ${\cal A}$-friendly labelings.  A labeling $f:V(G)\to {\cal A}$ is said to be {\em ${\cal A}$-friendly } if given any $a,b\in{\cal A}$, $-1\leq |f^{-1}(b)|-|f^{-1}(a)|\leq 1$.  If $g$ is an induced edge labeling and $f$ and $g$ are both ${\cal A}$-friendly Then $g$ is said to be an ${\cal A}$-cordial labeling and $G$ is said to be ${\cal A}$-cordial.  When ${\cal A}= \zk$ we say that $G$ is $k$-cordial.  We shall use this concept with digraphs.

Let $D=(V,A)$ be a  directed graph with vertex set $V$ and arc set $A$.   Let $f:V\to\{0,1\}$ be a friendly labeling of the vertices of $D$.  As for undirected graphs, an  induced labeling of the arc set  is a mapping $g:A\to {\cal X}$ for some set ${\cal X}$ where for an arc $(u,v)=\overrightarrow{uv}, g(u,v)= \hat{g}(f(u),f(v))$ for some $\hat{g}:\{0,1\}\times\{0,1\}\to {\cal X}$.  As we are dealing with directed graphs, it would be desirable for the induced labeling to distinguish between the label of the arc $(u,v)$  and the label of the arc $(v,u)$,  otherwise, the labeling would be an induced labeling of the underlying undirected graph.    If we let ${\cal X}={\mathbb Z}_3=\{-1,0,1\}$ and $\hat{g}(f(u),f(v))=f(v)-f(u)$,  we have an asymmetric labeling.  In this case, if   the set of arcs are nearly equally distributed among the three labelings,  we say that the labeling is $(2,3)$-cordial.   Formally:

\D  Let $D\in\Tn$, $D=(V,A)$.  Let  $f:V\to \{0,1\}$ be a friendly labeling of the vertex set $V$ of $D$.  Let $g:A\to\{1,0,-1\}$ be an induced labeling of the arcs of $D$.   If $g$ is friendly, that is, for any $i,j\in\{1,0,-1\}$, $-1\leq |g^{-1}(i)|-|g^{-1}(j)| \leq 1$.  such a labeling is called a \underline{\em $(2,3)$-cordial labeling}, and a digraph $D\in\Tn$ that can possess a $(2,3)$-cordial labeling will be called a \underline{\em$(2,3)$-cordial digraph.}
\DD

For any digraph $D$ with $(0,1)$-vertex labeling $f$, and induced  $(0,1,-1)$-arc labeling $g$, let $\Lambda_{f,g}(D)=(\alpha,\beta,\gamma)$ where $\alpha=|g^{-1}(1)|, \beta=|g^{-1}(-1)|,$ and $\gamma=|g^{-1}(0)|$.

Let $D\in\Tn$ and let $D^R$ be the digraph such that every arc of $D$ is reversed, so that $\overrightarrow{uv}$ is an arc in $D^R$ if and only if $\overrightarrow{vu}$ is an arc in $D$.  Let $f$ be a $(0,1)$-labeling of the vertices of $D$ and let $g(\overrightarrow{uv}) =f(v)-f(u)$ so that $g$ is a $(1,-1,0)$-labeling of the arcs of $D$.  Let $\overline{f}$ be the complementary $(0,1)$-labeling of the vertices of $D$, so that $\overline{f}(v)=0$ if and only if $f(v)=1$.  Let $\overline{g}$ be the corresponding induced arc labeling of $D$, $\overline{g}(\overrightarrow{uv}) =\overline{f}(v)-\overline{f}(u)$.

\l \label{lab} Let $D\in\Tn$ with vertex labeling $f$ and induced arc labeling $g$.  Let $\Lambda_{f,g}(D)=(\alpha,\beta,\gamma)$.   Then \begin{enumerate}\item $\Lambda_{f,g}(D^R)=(\beta,\alpha,\gamma)$. \item $\Lambda_{\overline{f},\overline{g}}(D)=(\beta,\alpha,\gamma)$, and \item $\Lambda_{\overline{f},\overline{g}}(D^R)=\Lambda_{f,g}(D).$\end{enumerate} \ll
\d If an arc is labeled 1, -1, 0 respectively then reversing the labeling of the incident vertices gives a labeling of -1, 1, 0 respectively,  If an arc $\overrightarrow{uv}$ is labeled 1, -1, 0 respectively, then $\overrightarrow{vu}$ would be labeled -1, 1, 0 respectively. \dd

\section{(2,3)-Cordial Digraphs.}

In this section, we shall use the arc labeling function $g:A\to \{0,1,-1\}$ defined by $g(\overrightarrow{v_iv_j})=g(v_i,v_j)=f(v_j)-f(v_i)$ where $f$ is a $(0,1)$-labeling of the vertices.

\subsection{(2,3)-cordial labelings of orientations of small trees}

We now give labelings of orientations of small trees to be used in inductive proofs later.

Let $P_n$ be a path of order $n$ (the number of vertices, so that there are $n-1$ arcs in $P_n$).  For $n=1$, there is only one directed graph and it is vacuously (2,3)-cordial.  For $n=2$ the only tree is an edge graph and any orientation is the digraph of an arc, which is (2,3)-cordial if the labelings of the two vertices are different.  

For $n=3$, every tree is a 2-path and so there are only two (non-isomorphic) orientations.  We list them here with (2,3)-cordial labelings:

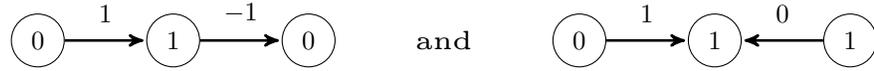
\begin{figure}[h]\label{3path}
\begin{center}
\begin{tikzpicture}[scale=0.9]

\tikzset{vertex/.style = {shape=circle,draw,minimum size=2em}}
\tikzset{edge/.style = {->,> = stealth',shorten >=1pt,thick}}

\node[vertex] (v1) at  (0,6) {$0$};
\node[vertex] (v2) at  (2,6) {$1$};
\node[vertex] (v3) at  (4,6) {$0$};
\path (6,6) node  [scale=2]     (x) {
\tiny
 and} ;
\node[vertex] (v2k-2) at  (8,6) {$0$};
\node[vertex] (v2k-1) at  (10,6) {$1$};
\node[vertex] (v2k) at  (12,6) {$1$};

\draw[edge,->, line width=1.0pt] (v1) to (v2);
\draw[edge,->, line width=1.0pt] (v2) to (v3);
\draw[edge,->, line width=1.0pt] (v2k-2) to (v2k-1);
\draw[edge,->, line width=1.0pt] (v2k) to (v2k-1);

\path (1,6.4) node     (y1) {$1$};
\path (3,6.4) node     (y2) {$-1$};
\path (9,6.4) node     (y5) {$1$};
\path (11,6.4) node     (y6) {$0$};

\end{tikzpicture}

\end{center}

  \caption{ $(2,3)$-Cordial labelings of orientations of 3-paths }

\end{figure}

For $n=4$ there are two undirected trees with three different (non-isomorphic) orientations of the 4-path and two non-isomorphic orientations of the star on 4 vertices.  Here we give these five oriented trees with (2,3)-cordial labelings of three of them:  See Figures 2 and 3.

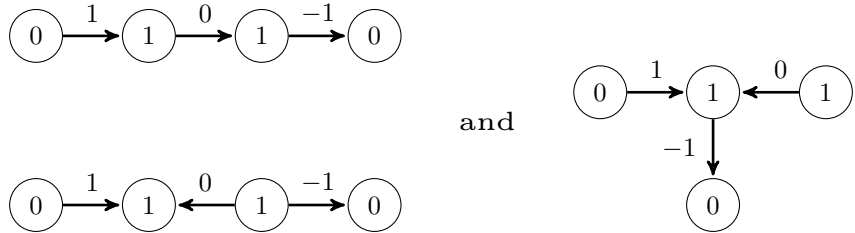
\begin{figure}[h]\label{fig3}
\begin{center}
\begin{tikzpicture}[scale=0.75]

\tikzset{vertex/.style = {shape=circle,draw,minimum size=2em}}
\tikzset{edge/.style = {->,> = stealth',shorten >=1pt,thick}}

\node[vertex] (v1) at  (0,6) {$0$};
\node[vertex] (v2) at  (2,6) {$1$};
\node[vertex] (v3) at  (4,6) {$1$};
\node[vertex] (v4) at  (6,6) {$0$};
\node[vertex] (u1) at  (0,3) {$0$};
\node[vertex] (u2) at  (2,3) {$1$};
\node[vertex] (u3) at  (4,3) {$1$};
\node[vertex] (u4) at  (6,3) {$0$};
\path (8,4.5) node  [scale=2]     (x) {
\tiny
 and} ;
\node[vertex] (v2k-3) at  (10,5) {$0$};
\node[vertex] (v2k-2) at  (12,5) {$1$};
\node[vertex] (v2k-1) at  (14,5) {$1$};
\node[vertex] (v2k) at  (12,3) {$0$};

\draw[edge,->, line width=1.0pt] (v1) to (v2);
\draw[edge,->, line width=1.0pt] (v2) to (v3);
\draw[edge,->, line width=1.0pt] (v3) to (v4);
\draw[edge,->, line width=1.0pt] (v2k-3) to (v2k-2);
\draw[edge,->, line width=1.0pt] (v2k-1) to (v2k-2);
\draw[edge,->, line width=1.0pt] (v2k-2) to (v2k);
\draw[edge,->, line width=1.0pt] (u1) to (u2);
\draw[edge,->, line width=1.0pt] (u3) to (u2);
\draw[edge,->, line width=1.0pt] (u3) to (u4);

\path (1,6.4) node     (y1) {$1$};
\path (3,6.4) node     (y2) {$0$};
\path (5,6.4) node     (y3) {$-1$};
\path (11,5.4) node     (y4) {$1$};
\path (13.2,5.4) node     (y5) {$0$};
\path (11.4,4) node     (y6) {$-1$};
\path (1,3.4) node     (y1) {$1$};
\path (3,3.4) node     (y2) {$0$};
\path (5,3.4) node     (y3) {$-1$};

\end{tikzpicture}

\end{center}

  \caption{ $(2,3)$-Cordial labeling of three 4-trees }

\end{figure}

\vskip 1cm

\begin{figure}[h]\label{fig3}
\begin{center}
\begin{tikzpicture}[scale=0.75]

\tikzset{vertex/.style = {shape=circle,draw,minimum size=2em}}
\tikzset{edge/.style = {->,> = stealth',shorten >=1pt,thick}}

\node[vertex] (v1) at  (0,6) {};
\node[vertex] (v2) at  (2,6) {};
\node[vertex] (v3) at  (4,6) {};
\node[vertex] (v4) at  (6,6) {};
\path (8,6) node  [scale=2]     (x) {
\tiny
 and} ;
\node[vertex] (v2k-3) at  (10,6) {};
\node[vertex] (v2k-2) at  (12,6) {};
\node[vertex] (v2k-1) at  (14,6) {};
\node[vertex] (v2k) at  (12,4) {};

\draw[edge,->, line width=1.0pt] (v1) to (v2);
\draw[edge,->, line width=1.0pt] (v2) to (v3);
\draw[edge,->, line width=1.0pt] (v4) to (v3);
\draw[edge,->, line width=1.0pt] (v2k-3) to (v2k-2);
\draw[edge,->, line width=1.0pt] (v2k-1) to (v2k-2);
\draw[edge,->, line width=1.0pt] (v2k) to (v2k-2);

\end{tikzpicture}

\end{center}

  \caption{ Two 4-trees that are not $(2,3)$-cordial }

\end{figure}
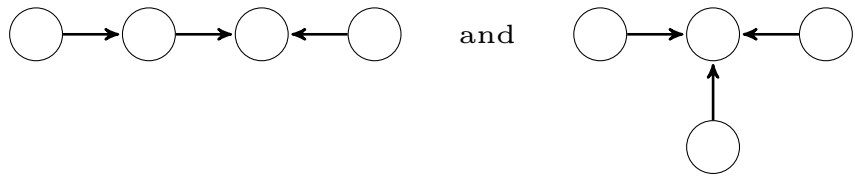

\subsubsection{(2,3)-Cordial Stars}  

Recall that for a friendly labeling $f$ of the vertices we let $g$ be the induced labeling of the arcs defined by $g(u,v)=f(v)-f(u)$.  

Let $S\in\Gn$ be a star graph with central vertex $v$ and $n-1$ pendant vertices.  A natural question is:  Is there an orientation of $S$ that is (2,3)-cordial?  Another is:  Are all orientations of $S$ (2,3)-cordial?  The answer to the latter question is:  Only if $n\leq 3$.    The answer to the first question is more complicated, but is completely answered here.

Let $n\geq 4$.  There are six cases to consider:  $n$ is even or odd and $n$ is a multiple of three, one more than a multiple of three or two more than a multiple of three.  Investigating these six cases, we let $k$ and $\ell$ be positive integers.

\noindent
{\bf Case 1.  $n$ is even, $n=2k$.}  Here there must be $k-1$ pendant vertices labeled the same as the labeling of the central vertex $v$.  Since $n\geq 4$, $k\geq 2$, we have that $\lceil\frac{k}{2}\rceil\leq k-1$.  For $g$ to be friendly, we must have $k-1$ arcs labeled 0, $\lceil\frac{k}{2}\rceil$ arcs labeled either 1 or -1 respectively, and $\lfloor\frac{k}{2}\rfloor$ arcs labeled -1 or 1 respectively.

\noindent
Subcase a.   The number of arcs is $n-1=3\ell$ so that $n=3\ell+1$.  For $S$ to be (2,3)-cordial we must have $k-1=\ell$.  But then $2k=n=3\ell+1 =3(k-1)+1$ so that $k=2$ and $n=4$.  

\noindent
Subcase b.   The number of arcs is $n-1=3\ell+1$, so that  $n=3\ell+2$.  For $S$ to be (2,3)-cordial we must have $k-1=\ell+1$ or   $k=\ell+2$. Hence $n=2k=2\ell+4=3\ell+2$, so that, $\ell=2$ , $k=4$ and $n=8$. 

\noindent
Subcase c.   The number of arcs is $n-1=3\ell +2$ so that  $n=3\ell+3$.  For $S$ to be (2,3)-cordial we must have $k-1=\ell+1$ or $k=\ell+2$.  Thus, $n=2k=2\ell+4=3\ell+3$ so that $\ell=1$ and $n=6$.

Thus, since any orientation of a 2-star is $(2,3)$-cordial, when $n$ is even we have that $n=2,4,6$ or $8$ only.

\noindent
{\bf Case 2. $n$ is odd, $n=2k+1$.}  Since for $n=1$ or 3, every orientation of a 1-star or 3-star is $(2,3)$-cordial, we assume that $n\geq 5$   and thus  $\lceil\frac{k}{2}\rceil\leq k-1$.   Hence there are  $k-1$ pendant vertices labeled the same as the labeling of the central vertex $v$.  So there are  $k-1$ arcs labeled 0.  

\noindent
Subcase a.  The number of arcs is  $n-1=3\ell$ so that  $n=3\ell+1$.  For $S$ to be (2,3)-cordial we must have  $k-1=\ell$ so that $k=\ell+1$ and hence $3\ell+1 = n = 2k+1=2\ell+3$ and hence $\ell=2$ and $n=7$.

\noindent
Subcase b.  The number of arcs is $n-1=3\ell+1$ so that $n=3\ell+2$.  For $S$ to be (2,3)-cordial we must have $k-1=\ell$ or $k=\ell +1$ and $\frac{k}{2}= \ell$ so that $n=5$; or $k-1=\ell+1$ and $k=\ell+2$, so that $n=11$.

\noindent
Subcase c.  The number of arcs is $n-1=3\ell +2$ so that $n=3\ell+3$.   For $S$ to be (2,3)-cordial we must have $k-1=\ell+1$ so $k=\ell+2$.  Thus, $3\ell+3=n=2k+1=2\ell+5$.  Hence $l=2$ and $n=9$.

Thus, when $n$ is odd, we have that $n=1,3,5,7,9$ or 11 only.

The above analysis  establishes:
\l
A star graph has a (2,3)-cordial orientation if and only if $n\leq 11$ and $n\neq 10$.  \ll

\subsubsection{(2,3)-Cordiality of Directed Cycles and Paths}

\begin{figure}[h]\label{fig4}
\begin{center}
\begin{tikzpicture}[scale=0.85]

\tikzset{vertex/.style = {shape=circle,draw,minimum size=3em}}
\tikzset{edge/.style = {->,> = stealth',shorten >=1pt,thick}}

\node[vertex] (v1) at  (0,6) {$v_1$};
\node[vertex] (v2) at  (2,6) {$v_2$};
\node[vertex] (v3) at  (4,6) {$v_3$};
\path (6,6) node  [scale=2]     (x) {$\dots$} ;
\node[vertex] (v2k-2) at  (8,6) {$v_{2k-2}$};
\node[vertex] (v2k-1) at  (10,6) {$v_{2k-1}$};
\node[vertex] (v2k) at  (12,6) {$v_{2k}$};

\path (0,7) node     (x1) {$0$};
\path (2,7) node     (x2) {$1$};
\path (4,7) node     (x3) {$0$};
\path (8,7) node     (x4) {$1$};
\path (10,7) node     (x5) {$0$};
\path (12,7) node     (x6) {$1$};

\path (1,5.4) node     (y1) {$1$};
\path (3,5.4) node     (y2) {$-1$};
\path (5,5.4) node     (y3) {$1$};
\path (7,5.4) node     (y4) {$1$};
\path (9,5.4) node     (y5) {$-1$};
\path (11,5.4) node     (y6) {$1$};

\draw[edge,->, line width=1.0pt] (v1) to (v2);
\draw[edge,->, line width=1.0pt] (v2) to (v3);
\draw[edge,->, line width=1.0pt] (v3) to (x);
\draw[edge,->, line width=1.0pt] (x) to (v2k-2);
\draw[edge,->, line width=1.0pt] (v2k-2) to (v2k-1);
\draw[edge,->, line width=1.0pt] (v2k-1) to (v2k);

\node[vertex] (v3k) at  (0,3) {$v_{3k}$};
\node[vertex] (v3k-1) at  (2,3) {$v_{3k-1}$};
\path (4,3) node  [scale=2]     (z) {$\dots$} ;

\node[vertex] (v2k+l+1) at  (6,3) {$ v_{5\ell+1}$};
\path (10,3) node  [scale=2]     (r) {$\dots$} ;
\node[vertex] (v2k+l) at  (8,3) {$v_{2k+\ell}$};
\node[vertex] (v2k+1) at  (12,3) {$v_{2k+1}$};

\draw[edge,->, line width=1.0pt] (v2k) to (v2k+1);
\draw[edge,->, line width=1.0pt] (v2k+1) to (r);
\draw[edge,->, line width=1.0pt] (r) to (v2k+l);
\draw[edge,->, line width=1.0pt] (v2k+l) to (v2k+l+1);
\draw[edge,->, line width=1.0pt] (v2k+l+1) to (z);
\draw[edge,->, line width=1.0pt] (z) to (v3k-1);
\draw[edge,->, line width=1.0pt] (v3k-1) to (v3k);
\draw[edge,->, line width=1.0pt] (v3k) to (v1);

\path (0,2) node     (Z1) {$0$};
\path (2,2) node     (Z2) {$0$};
\path (6,2) node     (Z3) {$0$};
\path (8,2) node     (Z4) {$1$};
\path (12,2) node     (x6) {$1$};

\path (1,3.6) node     (y1) {$0$};
\path (3,3.6) node     (y2) {$0$};
\path (5,3.6) node     (y3) {$0$};
\path (7,3.6) node     (y4) {$-1$};
\path (9,3.6) node     (y5) {$0$};

\path (.5,4.5) node     (y6) {$0$};
\path (11.5,4.5) node     (y6) {$0$};

\path (6,4.5) node     (y6) {\framebox[1.1\width]{$n=6\ell = 3k,  k=2\ell$}};

\end{tikzpicture}

\end{center}

  \caption{ $(2,3)$-Cordial labeling of a cycle }\label{cyc}

\end{figure}
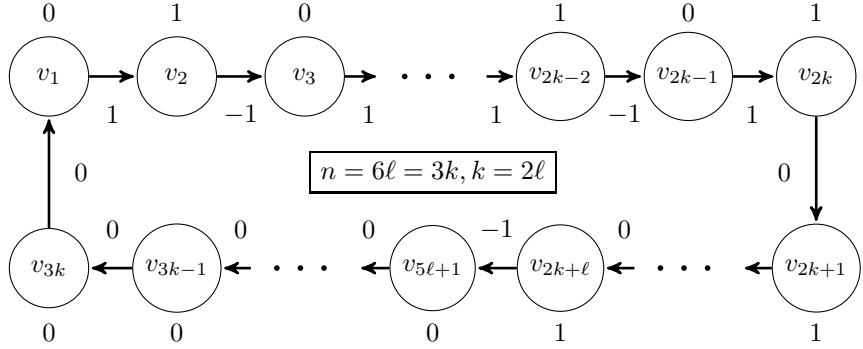

\t Let $f$ be a $(0,1)$-friendly labeling of $D\in\Tn$ and $g$ the induced $(1,0,-1)$-labeling of the arcs of $D$.  Then, 
\begin{enumerate}
\item  If $D$ is a directed $k$-cycle for $k\geq 3$, then $D$ is (2,3)-cordial.
\item If $D$ is any directed path in $\Dn$, then $D$ is (2,3)-cordial.  
\item  If $n\geq 4$ and $D$ is an out-star or an in-star then $D$ is not (2,3)-cordial.
\end{enumerate}
\tt
\d
The labeling of an $n$-cycle is shown in Figure \ref{cyc}, where $k=2\ell$ and $n=3k=6\ell$.  The vertex labeling is outside the cycle and the induced edge labeling is inside the cycle.  By deleting $1, 2, \dots, 5$ vertices, carefully chosen, one obtains a $(2,3)$-cordial labeling of a cycle of any length.  By deleting one edge, carefully chosen, one obtains a $(2,3)$-cordial labeling of any directed path.  Thus, conclusions 1 and 2 are proven.  

If $D$ is an  out-star with central vertex $v$, then, if the vertex $v$ is labeled "0", all arcs are labeled either "0" or "-1", while if the vertex $v$ is labeled "1", all arcs are labeled either "0" or "1"  Thus,  no labeling is (2,3)-cordial.  A parallel argument holds for an in-star. \dd

It is easily shown, even though time consuming, that both of the non-isomorphic tournaments on 3 vertices are (2,3)-cordial while only three of the four non-isomorphic tournaments on four vertices are cordial;  the 4-tournament with score vector (1,1,1,3) is not (2,3)-cordial.

\cj Given any tree $T \in \Gn$ with $\Delta(T)\leq 3$ there is an orientation of $T$ that is $(2,3)$-cordial. \cjj

\section{Paths}

\l\label{1-7} For $n\leq 9$, every orientation of any path on n vertices is (2,3)-cordial, except when $n=4$ and the orientation is one of those listed in Figure  5. \ll
\d  For $n=1$ or $n=2$ the lemma is trivially established.  For $n=3$ the two orientations in Figure 1 (and consequently their reversals) are shown to be (2,3)-cordial. 

For $n=4$, the four orientations of the 4-path that are not in Figure 5 are shown to be (2,3)-cordial in Figure 6.

For $n=5,6,$ and $7$, see the appendices,  There is a (2,3)-cordial labeling of each orientation of the paths of length 5, 6, and 7.

For $n=8$, append a vertex to the end of each orientation of the path on 7 vertices and obtain two orientations of the 8-path, one with $\overrightarrow{v_7v_8}$ and the other with  $\overleftarrow{v_7v_8}$.  By labeling vertex $v_8$ with 0 if there are fewer vertices in the labeling of the orientation of $P_7$ labeled 0 than 1 and labeling vertex $v_8$  a 1 otherwise we get a $(2,3)$-cordial labeling.  Thus,  we get a (2,3)cordial labeling of all 128 orientations of $P_8$.

For $n=9$, append a vertex to the end of each orientation of the path on 8 vertices and obtain two orientations of the 9-path, one with $\overrightarrow{v_8v_9}$ and the other with  $\overleftarrow{v_8v_9}$.  By choosing appropriately the label of $v_9$ one gets a (2,3)-cordial labeling.  This  gives a (2,3)-cordial labeling of all 256 orientations of $P_9$. \dd

\begin{figure}[h]
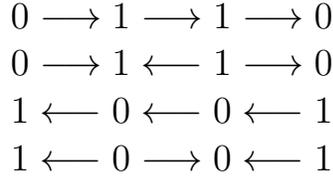
\label{4Cpaths}
\begin{center}
{\Large 

$0\>1\>1\>0$

$0\>1\<1\>0$

$1\<0\<0\<1$

$1\<0\>0\<1$

}

\end{center}

  \caption{ Labeled orientations of the 4-path that are  $(2,3)$-cordial }

\end{figure}

\begin{figure}[h]
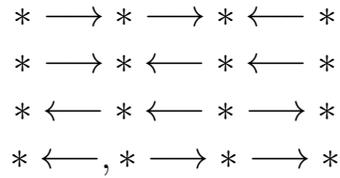
\label{4NCpaths}
\begin{center}

{\Large
$*\>*\>*\<*$

$*\>*\<*\<*$

$*\<*\<*\>*$

$*\<,*\>*\>*$

}

\end{center}

 \caption{ Orientations of the 4-path that are not $(2,3)$-cordial }

\end{figure}

Note that in the sequel, only vertex labelings are given as the arc labels are implicit.  For example in Figures \ref{4Cpaths} and \ref{4NCpaths} we give the 4-paths that are (2,3)-cordial  and the  4-paths that are not (2,3)-cordial.

\pagebreak

\cj Let $P$ be a path in $\Gn$ of length $n\geq 5$.  Then every orientation of $P$ is (2,3)-cordial.\cjj

The author wishes to thank the referee for his valuable suggestions.

\pagebreak

\section{APPENDIX: Orientations of paths of lengths 5, 6 and  7.}
 
\begin{figure}[h] 
\begin{multicols}{2}
\begin{center}
{\large
$0\>1\>1\>0\>0$

$0\>1\<1\>0\>0$

$0\>1\>0\<0\>1$

$0\>1\>1\>0\<0$

$0\>1\<1\<1\>0$

$0\>1\<1\>0\<0$

$0\>1\>1\<0\<1$

$0\>1\<0\<0\<1$}

\end{center}

\begin{center}
{\large
$1\<0\>0\>1\>0$

$1\<0\<0\>1\>0$

$1\<0\>0\<1\>1$

$1\<0\>0\>0\<1$

$1\<0\<0\<1\>1$

$1\<0\<0\>0\<1$

$1\<0\>0\<1\<1$

$1\<0\<0\<1\<1$}

\end{center}
\end{multicols}

 \caption{ Orientations of the 5-path  }
\end{figure}

\begin{figure}[h]
\begin{multicols}{2}
\begin{center}

$0\>1\>1\>0\>0\>1$

$0\>1\<1\>0\>0\>1$

$1\>0\<0\<1\>1\>0$

$1\<0\>0\<1\>1\>0$

$0\>1\>1\>0\<0\>1$

$0\>1\<1\>0\<0\>1$

$1\<0\<0\<1\<1\>0$

$1\<0\>0\<1\<1\>0$

$0\>1\>0\<0\>1\>1$

$0\>1\<0\<1\>1\>0$

$1\<0\<0\>1\>1\>0$

$1\<0\>1\>1\>0\>0$

$0\>0\>1\<0\<1\>1$

$0\>1\<0\<0\<1\>1$

$1\<1\<0\>0\<1\>0$

$1\<0\>0\>1\<1\>0$

\end{center}

\begin{center}
$0\>1\>1\>0\>0\<1$

$0\>1\<1\>0\>0\<1$

$1\>0\<0\<1\>1\<0$

$1\<0\>0\<1\>1\<0$

$0\>1\>1\>0\<0\<1$

$0\>1\<1\>0\<0\<1$

$1\<0\<0\<1\<1\<0$

$1\<0\>0\<1\<1\<0$

$0\>1\>0\<0\>1\<1$

$0\>1\<0\<1\>1\<1$

$1\<0\<0\>1\>1\<0$

$1\<0\>1\>1\>0\<1$

$0\>0\>1\<0\<1\<1$

$0\>1\<0\<0\<1\<1$

$1\<1\<0\>0\<1\<0$

$1\<0\>0\>1\<1\<0$
\end{center}
\end{multicols}
 \caption{ Orientations of the 6-path  }
\end{figure}

\begin{figure}[h]
\begin{multicols}{2}
\begin{center}

$0\rightarrow 1\rightarrow 1\rightarrow 0\rightarrow 0\rightarrow 1\rightarrow 0$

$0\rightarrow 1\leftarrow 1\rightarrow 0\rightarrow 0\rightarrow 1\rightarrow 0$

$0\leftarrow 1\leftarrow 1\leftarrow 0\rightarrow 0\rightarrow 1\rightarrow 0$

$0\leftarrow 1\rightarrow 1\leftarrow 0\rightarrow 0\rightarrow 1\rightarrow 0$

$0\rightarrow 1\rightarrow 1\rightarrow 0\leftarrow 0\rightarrow 1\rightarrow 0$

$0\rightarrow 1\leftarrow 1\rightarrow 0\leftarrow 0\rightarrow 1\rightarrow 0$

$0\leftarrow 1\leftarrow 1\leftarrow 0\leftarrow 0\rightarrow 1\rightarrow 0$

$0\leftarrow 1\rightarrow 1\leftarrow 0\leftarrow 0\rightarrow 1\rightarrow 0$

$0\rightarrow 1\rightarrow 1\rightarrow 0\rightarrow 1\rightarrow 0\leftarrow 0$

$0\rightarrow 1\leftarrow 1\rightarrow 0\rightarrow 1\rightarrow 0\leftarrow 0$

$0\leftarrow 1\leftarrow 1\leftarrow 0\rightarrow1\rightarrow0\leftarrow 0$

$0\leftarrow 1\rightarrow1\leftarrow 0\rightarrow1\rightarrow0\leftarrow 0$

$0\rightarrow1\rightarrow1\rightarrow0\leftarrow 1\rightarrow1\leftarrow 0$

$0\rightarrow1\leftarrow 1\rightarrow0\leftarrow 1\rightarrow1\leftarrow 0$

$0\leftarrow 1\leftarrow 1\leftarrow 0\leftarrow 1\rightarrow1\leftarrow 0$

$0\leftarrow 1\rightarrow1\leftarrow 0\leftarrow 1\rightarrow1\leftarrow 0$

$0\rightarrow1\rightarrow0\leftarrow 0\rightarrow1\rightarrow1\rightarrow0$

$0\rightarrow1\leftarrow 0\leftarrow 1\rightarrow1\rightarrow0\rightarrow0$

$0\leftarrow 1\leftarrow 0\rightarrow0\rightarrow1\rightarrow1\rightarrow0$

$0\leftarrow 0\rightarrow1\rightarrow0\rightarrow1\rightarrow1\rightarrow0$

$0\rightarrow1\rightarrow1\leftarrow 0\leftarrow 1\rightarrow1\rightarrow0$

$0\rightarrow1\leftarrow 0\leftarrow 0\leftarrow 1\rightarrow1\rightarrow0$

$1\leftarrow 0\leftarrow 0\rightarrow0\leftarrow 1\rightarrow0\rightarrow1$

$1\leftarrow 0\rightarrow0\rightarrow0\leftarrow 1\rightarrow0\rightarrow1$

$0\rightarrow1\rightarrow0\leftarrow 1\rightarrow1\rightarrow1\leftarrow 0$

$0\rightarrow0\leftarrow 1\leftarrow 1\rightarrow0\rightarrow1\leftarrow 0$

$0\leftarrow 0\leftarrow 1\rightarrow1\rightarrow0\rightarrow1\leftarrow 0$

$0\leftarrow 1\rightarrow1\rightarrow0\rightarrow0\rightarrow1\leftarrow 0$

$0\rightarrow1\rightarrow1\leftarrow 0\leftarrow 1\rightarrow0\leftarrow 0$

$0\rightarrow1\leftarrow 0\leftarrow 1\leftarrow 1\rightarrow0\leftarrow 0$

$1\leftarrow 0\leftarrow 0\rightarrow1\leftarrow 1\rightarrow0\leftarrow 1$

$1\leftarrow 0\rightarrow0\rightarrow1\leftarrow 1\rightarrow0\leftarrow 1$

$0\rightarrow0\rightarrow1\rightarrow0\rightarrow1\leftarrow1\rightarrow0$

$0\rightarrow1\leftarrow1\rightarrow0\rightarrow1\leftarrow1\rightarrow0$

$0\leftarrow0\leftarrow1\leftarrow0\rightarrow1\leftarrow1\rightarrow0$

$0\leftarrow1\rightarrow1\leftarrow0\rightarrow1\leftarrow1\rightarrow0$

$0\rightarrow1\rightarrow1\rightarrow0\leftarrow1\leftarrow0\rightarrow0$

$0\rightarrow1\leftarrow1\rightarrow0\leftarrow1\leftarrow0\rightarrow0$

$0\leftarrow1\leftarrow0\leftarrow1\leftarrow1\leftarrow0\rightarrow0$

$0\leftarrow1\rightarrow1\leftarrow0\leftarrow1\leftarrow0\rightarrow0$

$0\rightarrow1\rightarrow1\rightarrow0\rightarrow0\leftarrow1\leftarrow0$

$0\rightarrow1\leftarrow1\rightarrow0\rightarrow0\leftarrow1\leftarrow0$

$0\leftarrow1\leftarrow1\leftarrow0\rightarrow0\leftarrow1\leftarrow0$

$0\leftarrow1\rightarrow1\leftarrow0\rightarrow0\leftarrow1\leftarrow0$

$0\rightarrow1\rightarrow1\rightarrow0\leftarrow0\leftarrow1\leftarrow0$

$0\rightarrow1\leftarrow1\rightarrow0\leftarrow0\leftarrow1\leftarrow0$

$0\leftarrow1\leftarrow1\leftarrow0\leftarrow0\leftarrow1\leftarrow0$

$0\leftarrow1\rightarrow1\leftarrow0\leftarrow0\leftarrow1\leftarrow0$

$0\rightarrow1\rightarrow0\leftarrow0\rightarrow1\leftarrow1\rightarrow0$

$0\rightarrow0\leftarrow1\leftarrow0\rightarrow1\leftarrow1\rightarrow0$

$0\leftarrow1\leftarrow1\rightarrow0\rightarrow1\leftarrow0\rightarrow0$

$0\leftarrow1\rightarrow1\rightarrow0\rightarrow1\leftarrow0\rightarrow0$

$0\rightarrow1\rightarrow1\leftarrow0\leftarrow0\leftarrow1\rightarrow0$

$0\rightarrow1\leftarrow0\leftarrow0\leftarrow1\leftarrow1\rightarrow0$

$0\leftarrow1\leftarrow0\rightarrow1\leftarrow1\leftarrow1\rightarrow0$

$0\leftarrow1\rightarrow0\rightarrow1\leftarrow1\leftarrow0\rightarrow0$

$0\rightarrow0\rightarrow1\leftarrow1\rightarrow0\leftarrow1\leftarrow0$

$0\rightarrow1\leftarrow1\leftarrow1\rightarrow0\leftarrow1\leftarrow0$

$0\leftarrow1\leftarrow1\rightarrow0\rightarrow1\leftarrow1\leftarrow0$

$0\leftarrow1\rightarrow1\rightarrow0\rightarrow1\leftarrow1\leftarrow0$

$0\rightarrow1\rightarrow0\leftarrow0\leftarrow1\leftarrow1\leftarrow0$

$0\rightarrow0\leftarrow1\leftarrow0\leftarrow1\leftarrow1\leftarrow0$

$0\leftarrow0\leftarrow1\rightarrow1\leftarrow0\leftarrow1\leftarrow0$

$0\leftarrow1\rightarrow0\rightarrow1\leftarrow1\leftarrow0\leftarrow0$

\end{center}
\end{multicols}
\caption{ Orientations of the 7-path  }
\end{figure}

\end{document}